\documentclass[a4paper,12pt]{article}
\usepackage{amssymb,amsmath,amsthm}
\usepackage{graphicx}
\usepackage{url}

\newtheorem{theorem}{Theorem}[section]

\newtheorem{remark}{Remark}[section]

\def\E{{\bf E}}
\def\V{{{\bf Var}}}
\def\de{{\rm d}}
\def\I{{{\bf 1}}}

\numberwithin{equation}{section}
\title{Estimation for the discretely observed telegraph process}
\author{Iacus S.M.\footnote{Department of Economics, Business and Statistics, University of Milan, Via Conservatorio 7, 20122 Milan, Italy}, Yoshida N.\footnote{Graduate School of Mathematical Sciences, University of Tokyo, 3-8-1 Komaba, Meguro-ku, Tokyo 153-8914 Japan}}

\begin{document}
\maketitle

\begin{abstract}
The telegraph process $\{X(t), t>0\}$, is supposed to be  observed at $n+1$ equidistant time points
$t_i=i\Delta_n,i=0,1,\ldots, n$. The unknown value of  $\lambda$, the underlying rate of the Poisson process, is a parameter to be estimated. The asymptotic framework considered is the following:
 $\Delta_n \to 0$, $n\Delta_n = T \to \infty$ as $n \to \infty$.
We show that previously proposed moment type estimators are consistent and asymptotically normal but not efficient. We study further an approximated moment type estimator which is still not efficient but comes in explicit form. For this estimator the additional assumption $n\Delta_n^3 \to 0$ is required in order to obtain asymptotic normality. Finally, we propose a new estimator which is consistent, asymptotically normal and asymptotically efficient under no additional hypotheses.
\\\
\noindent {\bf key words:} telegraph process, discretely observed
process, inference for stochastic processes.
\end{abstract}

\section{Introduction}
The telegraph process (see Goldstein, 1951 and Kac, 1974) models
a random motion with finite velocity and it is usually proposed as an alternative to diffusion models.
The process
 describes the position of a particle moving
on the real line, alternatively with constant velocity $+ v$ or
$-v$. The changes of direction are governed by an homogeneous
Poisson process with rate $\lambda
>0.$ The telegraph process or  {\it telegrapher's} process is defined as
\begin{equation}\label{1.1}
X(t)=V(0)\int_0^t (-1)^{N(s)}\de s,\quad t>0,
\end{equation}
where  $V(0)$ is the initial velocity taking values  $\pm v$ with equal probability and independently of the Poisson process  $\{N(t), t>0\}$.
Many authors analyzed probabilistic properties of the process over the years (see
for example Orsingher, 1985, 1990; Pinsky, 1991; Foong and Kanno, 1994; Stadje and Zacks, 2004).
Di Crescenzo and Pellerey (2002) proposed the geometric telegraph
process as a model to describe the dynamics of the price of risky
assets, i.e. $S(t)=s_0\exp\{\alpha t+\sigma X(t)\}$, $t>0.$
where $X(t)$ replaces the standard Brownian motion of the original Black-Scholes (1973) - Merton (1973) model. Conversely to the standard geometric Brownian motion,
given that $X(t)$ is of bounded variation, so is the geometric telegraph process. This seems a realistic way to model paths of assets in the financial markets.
Mazza and Rulliere (2004) linked the process \eqref{1.1} and the
ruin processes in the context of risk theory. Di Masi {\it et al} (1994) proposed to model the volatility of financial markets in terms of the telegraph process. Ratanov (2004, 2005) proposed to model financial markets using a telegraph process with two intensities $\lambda_\pm$ and two velocities $c_\pm$. The telegraph process has also been used in  ecology to model population dynamics (see Holmes {\it et al.}, 1994) and the displacement of wild animals on the soil. In particular, this model is chosen because it preserves the property of animals to move at finite velocity and for a certain period along some direction (see e.g. Holmes, 1993, for an account).

It is worth to mention that, up to now, only few references about estimation problems for the telegrapher's  processes are  known.  Yao (1985) considers a the problem of state estimation of the telegrapher's process under white noise perturbation and studies performance of nonlinear filters.
Iacus (2001) is about the estimation of the parameter $\theta$ of the non-constant rate $\lambda_\theta(t)$ from continuous observations of the process.
More recently, De Gregorio and Iacus (2006) proposed pseudo-maximum likelihood and moment based estimators for the telegraph process under discrete observations on a fixed time interval $[0,T]$ when the process is observed with a mesh descreasing to zero.

The aim of this paper is the estimation of the parameter $\lambda$
when $\{X(t), 0 \leq t \leq T\}$ is observed at
equidistant times $0=t_0<\ldots<t_n$. We assume that
$t_i=i \Delta_n$, $i=0, \ldots, n$, hence $n\Delta_n=T$. The
asymptotic framework is the following: $\Delta_n \rightarrow 0$ and $n\Delta_n=T\to \infty$ as   $n\to \infty$.

When the telegraph process $X(t)$ is observed continuously then $N(T)/T$ is the optimal estimator of the parameter $\lambda$ and the statistical experiment is equivalent to the one of the observation of the whole  Poisson process on $[0,T]$ (see e.g.  Kutoyants, 1998). This situation also corresponds to the limiting  experiment in our asymptotic framework.

The paper is organized as follows. Section \ref{sec2} reviews some results on the telegraph process and presents a  formula of the $p$-th moment of the process in explicit form. This result is interesting in itself because it gives new information about this model.
Section \ref{sec3} presents  estimators previously introduced in the literature such as pseudo maximum likelihood estimators and moment type estimators. In particular,  for the moment type estimator 
it is shown that it is consistent and asymptotically normal but not efficient. Another approximated moment type estimator is given in explicit form and it is shown that the estimator is consistent and (under the additional assumption $n\Delta_n^3\to0$) asymptotically gaussian but still not efficient.
Finally, Section \ref{sec4} presents a new estimator which is consistent, asymptotically gaussian and asymptotically efficient without additional assumptions.

\section{Moments of the telegraph process}\label{sec2}
The process $X(t)$ is not Markovian. Conversely, the
two dimensional process $(X(t), V(t))$, $V(t) = V(0)(-1)^{N(t)}$, has the Markov property but
 a scheme of observation in which one is able
to observe both the position and the velocity of the process at
discrete time instants is not admissible, so statistical procedures should rely only on the observation of
the $X(t)$ component. The telegraph process is such that
\begin{equation}
\E X(t) = 0
\label{eq:m1}
\end{equation}
and
\begin{equation}
\E X^2(t)
=\frac{v^2}{\lambda}\left(t-\frac{1-e^{-2\lambda
t}}{2\lambda}\right)
\label{eq:m2}
\end{equation}
(see Orsingher, 1990). Next theorem gives the general explicit derivation of the moments of any order of the telegraph process. In the cited reference the author mentions that derivation of the moments of any order can be obtained but the actual derivation was not presented. In some sense, next theorem completes Section 3 of Orsingher (1990) and we decided to present it here because it has some interest in itself.
\begin{theorem}
Let $p\geq 1$, then
\begin{equation}
\E X(t)^p =  (1 + (-1)^p)  (vt)^p \left(\frac{2}{\lambda t}\right)^\frac{p-1}{2} \Gamma\left(\frac{p+1}{2}\right)
\left\{
I_\frac{p+1}{2}(\lambda t) +I_\frac{p-1}{2}(\lambda t)
\right\} \frac{e^{-\lambda t}}{2} 
\label{eq:pmoment}
\end{equation}
\label{thm:mmt}
\end{theorem}
\begin{proof}
We start by rewriting the $p$-th moment of $X(t)$ as the sum of the two terms emerging from the discrete and the absolute continuous part of its density (see formula (27) in Orsingher (1990) or \eqref{2.2} below). Therefore we have the following representation
\begin{equation}
E X(t)^p = \int\limits_{-vt}^{+vt} x^p p(x, t) \de x + (1 + (-1)^p)   \frac{(v t)^p e^{-\lambda t}}{2}
\label{eq:mmt}
\end{equation}
The term $e^{-\lambda t}/2$ will also appear in the above integral, so we consider separately the two identities
\begin{equation}
\frac{\lambda}{v}\int_{-vt}^{+vt} x^p 
I_0\left(\frac{\lambda}{v}\sqrt{v^2 t^2 -x^2} \right)\de x =
(1+(-1)^p)(vt)^p\left(\frac{2}{\lambda t}\right)^{\frac{p-1}{2}}    \Gamma\left(\frac{p+1}{2}\right) I_{\frac{p+1}{2}}(\lambda t)
\label{eq:int1}
\end{equation}
and
\begin{equation}
\begin{aligned}
\lambda t \int_{-vt}^{+vt} x^p &\frac{I_1\left(\frac{\lambda}{v} \sqrt{v^2 t^2 -x^2}
\right)}{\sqrt{v^2 t^2 -x^2}}\de x=\\ 
& = (1+(-1)^p) (vt)^p\left\{  \left(\frac{2}{\lambda t}\right)^\frac{p-1}{2}\Gamma\left(\frac{p+1}{2}\right) I_{\frac{p-1}{2}}(\lambda t) - 1\right\}
\label{eq:int2}
\end{aligned}
\end{equation}
After multiplication by the factor  $e^{-\lambda t}/2$ both \eqref{eq:int1} and \eqref{eq:int2}, direct substitution  in \eqref{eq:mmt} gives the result of the theorem.
So we need to prove the above identities and we start with  formula \eqref{eq:int1}
$$
\begin{aligned}
\frac{\lambda}{v}\int_{-vt}^{+vt} x^p 
& I_0\left(\frac{\lambda}{v}\sqrt{v^2 t^2 -x^2} \right)\de x
=
\frac{\lambda}{v}\int_{-vt}^{+vt} x^p 
\sum_{k=0}^\infty \frac{1}{(k!)^2} \left(\frac{\lambda}{2v}\sqrt{v^2 t^2 -x^2}
\right)^{2k}\de x\\
&=\frac{\lambda}{v}\sum_{k=0}^\infty \frac{1}{(k!)^2}  \left(\frac{\lambda}{2v}\right)^{2k}
\int_{-vt}^{+vt} x^p 
 \left(v^2 t^2 -x^2\right)^{k}\de x\\
 &=\frac{\lambda}{v}\sum_{k=0}^\infty \frac{1}{(k!)^2}  (1+(-1)^p)(vt)^{p+2k+1} 
 \frac{Beta\left(k+1,\frac{p+1}{2}\right)}{2}\\
&= \frac{\lambda}{v}2^{\frac{p-1}{2}}(1+(-1)^p) \left(\frac{v^2 t^2}{\lambda t}\right)^\frac{p+1}{2}  \Gamma\left(\frac{1+p}{2}\right) I_{\frac{p+1}{2}}(\lambda t)\\
&= (1+(-1)^p) \left(\frac{2}{\lambda t}\right)^{\frac{p-1}{2}} (vt)^p  \Gamma\left(\frac{1+p}{2}\right) I_{\frac{p+1}{2}}(\lambda t)
\end{aligned}
$$
where $Beta(a,b) = \Gamma(a)\Gamma(b)/\Gamma(a+b)$, $\Gamma(x) = \int_0^\infty e^{-t}t^{x-1}\de t$, $n! = \Gamma(n+1)$ and
$$
I_\nu(x) = \sum_{k=0}^\infty \frac{1}{\Gamma(k+1+\nu)k!}\left(\frac{x}{2}\right)^{2k+\nu}
$$
In the above, the derivation of the equality
$$
\int_{-vt}^{+vt} x^p 
 \left(v^2 t^2 -x^2\right)^{k}\de x=
 (1+(-1)^p)(vt)^{p+2k+1} 
 \frac{Beta\left(k+1,\frac{p+1}{2}\right)}{2}
 $$
is incredibly lengthy but trivial therefore we omit it. 
We now calculate \eqref{eq:int2}
$$
\begin{aligned}
\lambda t \int_{-vt}^{+vt} &x^p \frac{I_1\left(\frac{\lambda}{v} \sqrt{v^2 t^2 -x^2}
\right)}{\sqrt{v^2 t^2 -x^2}} 
\de x =\lambda t
\int_{-vt}^{+vt} x^p 
\sum_{k=0}^\infty\frac{1}{\Gamma(k+2)k!} \left(\frac{\lambda}{2v}\right)^{2k+1}   \left( v^2 t^2 -x^2
\right)^k\de x \\
&=\lambda t
\sum_{k=0}^\infty\frac{1}{\Gamma(k+2)k!} \left(\frac{\lambda}{2v}\right)^{2k+1}   
 (1+(-1)^p)(vt)^{p+2k+1} 
 \frac{Beta\left(k+1,\frac{p+1}{2}\right)}{2}\\
&=(1+(-1)^p) (vt)^p \frac{\lambda t}{2}  \Gamma\left(\frac{p+1}{2}\right) \left\{\left(\frac{2}{\lambda t}\right)^\frac{p+1}{2}I_{\frac12 (p-1)}(\lambda t) - \frac{1+p}{\lambda t \Gamma\left(\frac{3+p}{2}\right)}\right\}\\
&=(1+(-1)^p) (vt)^p\left\{  \Gamma\left(\frac{p+1}{2}\right) \left(\frac{2}{\lambda t}\right)^\frac{p-1}{2}I_{\frac12 (p-1)}(\lambda t) - 1\right\}
\end{aligned}
$$
because $ \Gamma((p+1)/2) / \Gamma((3+p)/2) = 2/(1+p)$.
\end{proof}
\begin{remark}\label{rem3}
From  \eqref{eq:pmoment} it emerges that all odd-moments of the process are identically zero. 
Moreover,  for $x \to 0$, the modified Bessel functions admit the following expansion 
$$I_\nu(x) =  \frac{1}{\Gamma(\nu+1)}\left(\frac{x}{2}\right)^\nu \left(
1+\frac{z^2}{4(\nu+1)} + \frac{z^4}{32(\nu+1)(\nu+2)}+ \cdots\right)
$$ 
from which we obtain that $E X(t)^p$ is of order $t^{p+2}$ for $t \to 0$.
The following expansion, for $t\to 0$, will be useful in the following
\begin{eqnarray}
\E X(t)^2  &=& v^2 t^2 -\frac23 v^2 \lambda t^3 + \frac13 v^2 \lambda^2 t^4 + o(t^4)\label{m2}\\
\E X(t)^4 &=& v^4 t^4 -\frac45 v^4  \lambda t^5 + \frac25 v^4 \lambda^2 t^6 + o(t^6)\label{m4}\\
\E X(t)^6  &=& v^6 t^6 -\frac67 v^6  \lambda t^7 + \frac37 v^6 \lambda^2 t^8 + o(t^8)\label{m6}
\end{eqnarray}
\end{remark}
We now check that for $p=2$ we recover formula \eqref{eq:m2} which has been derived in two different ways in Orsingher (1990).
Indeed, for $p=2$ we have
$$
2 (vt)^2 \left(\frac{2}{\lambda t}\right)^\frac{1}{2} \frac{\pi}{2} \frac{e^{-\lambda t}}{2} 
\left\{
I_\frac{3}{2}(\lambda t) +I_\frac{1}{2}(\lambda t)
\right\}
$$
and noticing that
$$
I_\frac{1}{2}(x) = \sqrt{\frac{2}{\pi}}\frac{\sinh(x)}{\sqrt{x}},\qquad
I_\frac{3}{2}(x) = \sqrt{\frac{2}{\pi}}\frac{x \cosh(x)-\sinh(x)}{\sqrt{x^3}}
$$
direct substitution gives \eqref{eq:m2}.
\begin{remark}\label{rem2}
The formula of the fourth moment has also a relatively simple expression, so we present it here. To derive the result, it is useful to know that
$$
I_\frac{5}{2}(x) = \sqrt{\frac{2}{\pi}}\frac{(x^2+3) \sinh(x)-3x\cosh(x)}{\sqrt{x^5}}
$$
Thus,
$$
\begin{aligned}
E X(t)^4 &=  2 (vt)^4 \left(\frac{2}{\lambda t}\right)^\frac{3}{2} \frac{3}{4}\sqrt{\pi}
\left\{
I_\frac{5}{2}(\lambda t) +I_\frac{3}{2}(\lambda t)
\right\} \frac{e^{-\lambda t}}{2} \\
&= 3 \left(\frac{v}{\lambda}\right)^4 e^{-\lambda t}\,
    \left\{ \lambda t\,\left( \lambda t -3 \right) \,
       \cosh (\lambda t ) + 
      \left( 3 + \lambda t\,
          \left(  \lambda t  -1 \right)  \right) \,
       \sinh (\lambda t) \right\}     
\end{aligned}
$$
\end{remark}
\section{Previous results on the estimation of $\lambda$}\label{sec3}
As mentioned in the Introduction, we assume that the telegraph process $\{X(t), 0 \leq t \leq T\}$,
with $X(0)=x_0=0,$ is observed only at discrete times $0 <
t_1 < \cdots <t_n=T$, with $t_i = i\Delta_n$, $i=0,
\ldots, n$ hence $n\Delta_n = T$. We use the following notation to
simplify the formulas: $X(t_i) = X(i\Delta_n) = X_i$.
The interest is in
the estimation of the parameter $\lambda$ whilst $v$ is assumed to
be known.
If  the whole trajectory can be oserved, $\lambda$ can be
estimated by $N(T)/T$ where $N(T)$ 
is the number of Poisson events counted in $[0,T]$ or,
 the number of times the
process switches its velocity in $[0,T]$.
The estimation of $v$ is always an uninteresting problem as, if
there are no switchings in $((i-1)\Delta_n, i\Delta_n]$ then $X_i
- X_{i-1} = v \Delta_n$, hence if $\Delta_n$ is sufficiently
small, there is high probability of observing $N(t_{i+1}) - N(t_i)
= 0$ then $v$ can be estimated (actually calculated) without
error.
The asymptotic minimum variance of all the estimators for the continuous time experiment is the value of $\lambda$ itself because, as said, it is just the problem of estimating the intensity of a homogeneous Poisson process.
We now review some estimators for this process already available in the literature and study some new properties of one of them.
De Gregorio and Iacus (2006) considered the following approximated likelihood 
\begin{eqnarray}
L_n(\lambda)&=&L_n(\lambda | X_0, X_1, \ldots, X_n) =
\prod_{i=1}^n
p(X_i, \Delta_n;X_{i-1},t_{i-1})\label{2.2}\\
&=&\prod_{i=1}^n\Bigg\{ \frac{e^{-\lambda\Delta_n}}{2v} \left\{
\lambda I_0\left(\frac{\lambda}{v}\sqrt{u_{n,i}} \right) + \frac{v
\lambda \Delta_n I_1\left(\frac{\lambda}{v} \sqrt{u_{n,i}}
\right)}{\sqrt{u_{n,i}}}
\right\}\I_{\{u_{n,i}>0\}}\notag\\
&&+\frac{e^{-\lambda \Delta_n}}{2} \delta(u_{n,i}=0)\Bigg\}\notag
\end{eqnarray}
where $u_{n,i} =  u_n(X_i,X_{i-1}) = v^2 \Delta_n^2 -  (
X_i-X_{i-1} )^2$, $\delta$ is the Dirac function and $\I_A$ is the indicator function of set $A$.
The density $p(X_i, \Delta_n;X_{i-1},t_{i-1})$ appearing in
\eqref{2.2} is  the probability law
of a telegraph process initially located in $X_{i-1}$, that
reaches the position $X_{i}$ at time $t_i$. The above approximated likelihood is indeed the joint law of the  increments $X_i - X_{i-1}$ which are considered as if they were $n$ independent copies of the process $X(\Delta_n)$. The increments $\eta_i$ can be expressed as follows
$$
\eta_i = X_i - X_{i-1} =  V(0) \int_{t_{i-1}}^{t_i} (-1)^{N(s)} \de s =
V(0)(-1)^{N(t_{i-1})} \int_{t_{i-1}}^{t_i} (-1)^{N(s)-N(t_{i-1})} \de s   
$$
and they are stationary but not independent. Conversely, the squared increments
$$
\eta_i^2 = v^2 \left( \int_{t_{i-1}}^{t_i} (-1)^{N(s)-N(t_{i-1})} \de s \right)^2  
$$
(or the absolute increments $|\eta_i|$) are independent.
In their paper, the authors proposed the following  estimator
\begin{eqnarray}\label{2.4}
\hat\lambda_n = \arg\max_{\lambda>0} L_n(\lambda)
\end{eqnarray}
The estimator is proved to be unique and to exist (not so evident given the uncommon form of $L_n$) and such that $\bar\lambda_n\rightarrow N(T)/T$
under the condition $n\Delta_n=T$, $\Delta_n\rightarrow 0$ as $n\to 0$ but $T$ fixed. The limiting estimator $ N(T)/T$ is  the natural estimator,  but $\bar\lambda_n$ is not consistent for all values of $\lambda$ because time $T$ is fixed. 
In the same paper, the authors present numerical results about a least squares estimator of the following form
\begin{equation}
\check{\lambda}_n = \arg\min_{\lambda>0} \left\{\frac{1}{n} \sum_{i=1}^n \eta_i^2 - \frac{v^2}{\lambda}\left(\Delta_n-\frac{1-e^{-2\lambda
\Delta_n}}{2\lambda}\right)\right\}^2
\label{deg}
\end{equation}
In order to have consistency and hence asymptotic normality of estimators it is necessary to consider the asymptotics as $n\Delta_n=T\to \infty$.
We will prove that the estimator $\check{\lambda}_n$ is a true moment type estimator which is consistent and asymptotically gaussian but not efficient. The estimator $\check{\lambda}_n$ is given in implicit form and we also study an approximated moment type estimator which is given in explicit form and prove
that it is consistent and asymptotically normal (under the additional condition $n\Delta_n^3\to0$) but still not efficient because it asymptotic variance is $\frac65 \lambda$. A new asymptotically efficient estimator will be presented in Section \ref{sec4}.

\subsection{The moment type estimator}
Consider the original estimator from \eqref{deg}. The statistics
$$U_n = \frac{1}{n} \sum_{i=1}^n \eta_i^2$$
is an unbiased estimator of
$$
u_0 = f(\lambda_0) = \frac{v^2}{\lambda_0}\left(\Delta_n-\frac{1-e^{-2\lambda_0
\Delta_n}}{2\lambda_0}\right)
$$ 
Observe now that $f(\lambda)$ is monotonic and decreasing function of $\lambda$ such that
$$
\lim_{\lambda\to 0} f(\lambda) = v^2 \Delta_n^2,\quad 
\lim_{\lambda\to \infty} f(\lambda) =0
$$
on the other side  $U_n$ varies from $0$ to $v^2 \Delta_n^2$, hence the minimum value of  \eqref{deg} $\check\lambda_n$ is also the solution of
\begin{equation}
\frac{1}{n} \sum_{i=1}^n \eta_i^2 - \frac{v^2}{\lambda}\left(\Delta_n-\frac{1-e^{-2\lambda
\Delta_n}}{2\lambda} \right)= 0
\label{eq:mmt2}
\end{equation}
which means that $\check\lambda_n$ is a true moment type estimator.
Let $\lambda_0$ be the true value of the parameter and $\E_0$ and $\V_0$ indicate the expected value and the variance operator under the true parameter $\lambda_0$. 
\begin{theorem}
Let $\check\lambda_n$ the moment type estimator solution of \eqref{eq:mmt2} and let $n\Delta_n\to \infty$, $\Delta_n\to0$ as $n\to\infty$ and let $\lambda_0>0$ be the true value of the parameter. Then, $\check\lambda_n$ is a consistent estimator of $\lambda_0$ and such that
$$
\sqrt{n\Delta_n}(\check\lambda_n - \lambda_0) \overset{d}{\to} N\left(0, \frac{6}{5}\lambda_0\right)
$$
where $\overset{p}{\to}$ denotes the convergence in  distribution.
\end{theorem}
\begin{proof}
We  rewrite $U_n$ as $f(\check\lambda_n)$ and study the asymptotic properties of $\check\lambda_n$ by $\delta$-method.
Let $\lambda=g(u) = f^{-1}(u)$ with $f^{-1}$ the inverse function of $f(\lambda)$. Hence
$$
\sqrt{n\Delta_n}(\check\lambda_n - \lambda_0) =\sqrt{n\Delta_n}(f^{-1}(U_n) - f^{-1}(u_0)) = 
\sqrt{n\Delta_n}(U_n - u_0) \frac{1}{f'(\lambda_0)} + o_p\left(\sqrt{n\Delta_n}|U_n -u_0|\right) 
$$
where
$$
\begin{aligned}
f'(\lambda) &= \frac{\de}{\de\lambda}f(\lambda) = 
- \frac{v^2\,\left( 1 + \Delta_n \,\lambda  + 
        e^{2\,\Delta_n \,\lambda }\,
         \left( -1 + \Delta_n \,\lambda  \right)  \right) }
      {e^{2\,\Delta_n \,\lambda }\,{\lambda }^3} 
\end{aligned}
$$
Then $\E_0 \sqrt{n\Delta_n}(\check\lambda_n - \lambda_0) = 0$. Moreover,
$$
\begin{aligned}
\V_0(\eta_i^2) &=\frac{v^4\,\left( -1 - 16\,e^{2\,t\,\lambda_0 }\,
       \left( 1 + \Delta_n\,\lambda_0  \right)  + 
      e^{4\,\Delta_n\,\lambda_0 }\,
       \left( 17 + 4\,\Delta_n\,\lambda_0 \,
          \left( -5 + 2\,\Delta_n\,\lambda_0  \right)  \right) 
      \right) }{4\,e^{4\,\Delta_n\,\lambda_0 }\,{\lambda_0 }^4}\\
        &\simeq \frac{8}{15} v^4 \lambda_0 \Delta_n^5 - \frac{32}{45} v^4 \lambda_0^4 \Delta_n^6 + o(\Delta_n^7)
\end{aligned}
$$
$$
\begin{aligned}
\V_0 &(\sqrt{n\Delta_n}(\check\lambda_n - \lambda_0)) \\
&= 
n\Delta_n \frac{1}{\left(f'(\lambda_0)\right)^2}\frac{1}{n}\V_0(\eta_i^2)\\
&=
\frac{\Delta_n\,{\lambda_0 }^2\,\left( -1 - 
      16\,e^{2\,\Delta_n\,\lambda_0 }\,
       \left( 1 + \Delta_n\,\lambda_0  \right)  + 
      e^{4\,\Delta_n\,\lambda_0 }\,
       \left( 17 + 4\,\Delta_n\,\lambda_0 \,
          \left( -5 + 2\,\Delta_n\,\lambda_0  \right)  \right) 
      \right) }{4\,{\left( 1 + \Delta_n\,\lambda_0  + 
        e^{2\,\Delta_n\,\lambda_0 }\,
         \left( -1 + \Delta_n\,\lambda_0  \right)  \right) }^2}\\
&\simeq\frac{6\,\lambda_0 }{5} + \frac{4\,{\lambda_0 }^2\,\Delta_n}{5} + o(\Delta_n^2)\to \frac{6}{5}\lambda_0
\end{aligned}
$$
Let
$$
\xi_i = \frac{\sqrt{n\Delta_n}}{nf'(\lambda_0)} (\eta_i^2 - f(\lambda_0))
$$
in order to prove asymptotic normality, we need to prove the Lindeberg condition. Therefore, we need to prove that $n\E_0 |\xi|^3\to 0$. Indeed,
$$
\begin{aligned}
n\E_0 |\xi_i|^3 &= \frac{\Delta_n^\frac32}{\sqrt{n}(f'(\lambda_0))^3} \E_0 |\eta_i^2 - f(\lambda_0)|^3\\
&= \frac{\Delta_n^\frac32}{\sqrt{n}(f'(\lambda_0))^3} \left(|v^2\Delta_n^2 -  f(\lambda_0)|^3 + o(\Delta_n^2)\right)\\
&\simeq - \frac{\lambda_0^3 \Delta_n^\frac32}{\sqrt{n}} + o(\Delta_n^\frac52)\to 0
\end{aligned}
$$
\end{proof}
\subsection{An approximated (but explicit) moment type estimator}
Consider again \eqref{eq:m2}. Some algebra, or Remark \ref{rem3}, give the following expansion
$$
\begin{aligned}
\E\left(X_{i} - X_{i-1}\right)^2 &= \frac{v^2}{\lambda}\left(\Delta_n-\frac{1-e^{-2\lambda
\Delta_n}}{2\lambda}\right)\\
&= \frac{v^2}{\lambda}
\left(
\Delta_n-\frac{2\lambda \Delta_n - \frac12 (-2\lambda \Delta_n)^2) - \frac16 (-2\lambda \Delta_n)^3) + o(\Delta_n^3)}{2\lambda}
\right)
\\
&=v^2 \Delta_n^2 - \frac23 v^2 \lambda \Delta_n^3 + o(\Delta_n^3)
\end{aligned}
$$
Therefore, an approximated moment type estimator is the following
\begin{equation}
\lambda_n^* = \frac{3}{2}\frac{1}{n v^2\Delta_n^3}\sum_{i=1}^n\left\{v^2\Delta_n^2- (X_i - X_{i-1})^2\right\}
 = \frac{3}{2}\frac{1}{n\Delta_n}\sum_{i=1}^n\left\{1- \frac{\eta_i^2}{v^2\Delta_n^2}\right\}
\label{est:naka}
\end{equation}
and $\lambda_n^*$ is a weighted sum of the independent random variables $\eta_i^2$.
Remark that $v^2\Delta_n^2- (X_i - X_{i-1})^2$ is exactly zero if no Poisson event occurs in the time interval $(t_i, t_{i+1}]$. This fact will be used to evaluate expected values of related quantities in the following. This estimator is qualitatively not different from the estimator $\check{\lambda}_n$ in equation \eqref{deg} but the cost of having an estimator in explicit form is paid by the need of use of the additional condition $n\Delta_n^3\to0$ in order to obtain asymptotic normality.
This hypothesis has been also used in the ``high frequency'' sampling  for discretely observed diffusion processes (see e.g. Florens-Zmirou, 1989; Yoshida, 1992). 

\begin{theorem}\label{th:gauss}
Let $\lambda_0$ be the true value of the parameter. Then, under the condition $\Delta_n\to 0$,  $n\Delta_n \to \infty$ as $n\to \infty$ the statistics $\lambda_n^*$ in \eqref{est:naka} is consistent estimator of $\lambda_0$. Moreover, under the condition $n\Delta_n^3\to 0$ it is also asymptotically Gaussian, i.e.
$$
\sqrt{n\Delta_n}(\lambda_n^* - \lambda_0) \overset{d}{\rightarrow} N\left(0, \frac{6}{5} \lambda_0\right) 
$$
\end{theorem}
\begin{proof}
Consistency is trivial, indeed
$$
\E_0 \lambda^*_n = \lambda_0 + o(1)
$$
We now study its asymptotic distribution, in particular we are interested in the following random variable
$$
\sqrt{n\Delta_n}(\lambda_n^* - \lambda_0) = 
 \frac{3}{2 \sqrt{n}v^2\Delta_n^\frac52}\sum_{i=1}^n \left\{v^2 \Delta_n^2 - (X_i - X_{i-1})^2 - \frac23 \lambda_0 \Delta_n^3 \right\}
= \sum_{i=1}^n \xi_i\\
$$
We now show that
$$
\E_0 \left\{\sqrt{n\Delta_n}(\lambda_n^* - \lambda_0) \right\} \simeq \sqrt{n\Delta_n^3}
$$
holds true, where  $a \simeq b$ means ``$a$ is the same order of $b$''.
This is where the assumption $n\Delta_n^3\to0$ is needed in order to obtain asymptotic normality for this estimator.
$$
\begin{aligned}
\left|\E_0 \sum_i \xi_i\right| &\leq n |\xi_i|  \left(P\{N((t_{i-1},t_i])=0\} +  P\{N((t_{i-1},t_i])=1\}  +  P\{N((t_{i-1},t_i])>1\}\right) \\
&=\frac{3n}{2\sqrt{n}\Delta_n^\frac52} \left|-\frac23 \lambda_0\Delta_n^3 v^2 + \lambda_0 \Delta_n \left(
v^2 \Delta_n^2 -v^2 \frac{\Delta_n^2}{3} - \frac23 \lambda_0 \Delta_n^3\right) + o(\Delta_n^4)
\right|\\
&\simeq \frac{n\Delta_n^4}{\sqrt{n}\Delta_n^\frac52} = \sqrt{n\Delta_n^3}\to 0
\end{aligned}
$$
Moreover, ${\rm Cov}(\xi_i, \xi_j) = 0$, for $i\neq j$. 
Now we need to prove Lindeberg condition.
We  calculate the variance of $\xi_i$ making use of the following two expansions derived from \eqref{m2}, \eqref{m4} and \eqref{m6}
$$
\begin{aligned}
\E \eta_i^2 &= \E (X_i - X_{i-1})^2 = v^2 \Delta_n^2 -\frac23 \lambda \Delta_n^3 + \frac13 v^2 \lambda^2 \Delta_n^4 + o(\Delta_n^4)\\
\E \eta_i^4 &= v^4 \Delta_n^4 -\frac45 v^4  \lambda \Delta_n^5 + \frac25 v^4 \lambda^2 \Delta_n^6 + o(\Delta_n^6)\\
\E \eta_i^6  &= v^6 \Delta_n^6 -\frac67 v^6  \lambda \Delta_n^7 + \frac37 v^6 \lambda^2 \Delta_n^8 + o(\Delta_n^8)\\
\end{aligned}
$$
so that
$$
\V( \eta_i^2 )= \E_0 \eta_i^4 - \left(\E_0 \eta_i^2\right)^2 =
\frac{8}{15} v^4 \lambda_0 \Delta_n^5 - \frac{32}{45}v^4\lambda_0^2 \Delta_n^6 + o(\Delta_n^6)
$$
hence
$$
\V( \xi_i) = \frac{9}{4 v^4 n \Delta_n^5} \V (\eta_i^2) = \frac{6}{5}\frac{\lambda_0}{n}-
\frac{8}{5}  \lambda_0^2\frac{\Delta_n}{n} + \frac{9}{4n v^4} o(\Delta_n)
$$
Thus
\begin{equation}
\V\left\{ \sqrt{n\Delta_n}(\lambda_n^* - \lambda_0) \right\} = n \V (\xi_i) = 
 \frac{6}{5} \lambda_0 -
\frac{8}{5}  \lambda_0^2 \Delta_n + \frac{9}{4 v^4} o(\Delta_n)\to  \frac{6}{5} \lambda_0 
\label{eq:limvar}
\end{equation}
Now, let $S_n = \gamma_1+\gamma_2+\cdots+\gamma_n$ and $s^2_n = \V S_n$
with 
$$\gamma_i = \frac{3}{2 nv^2\Delta_n^3}\left\{v^2 \Delta_n^2 - (X_i - X_{i-1})^2 - \frac23 \lambda_0 \Delta_n^3 \right\} = \frac{\xi_i}{\sqrt{n\Delta_n}}
$$
We have that
$$s^2_n = \V S_n = n \V(\gamma_i) = \frac{1}{\Delta_n} \V (\xi_i)= 
 \frac{6}{5} \frac{\lambda}{n\Delta_n} -
\frac{8}{5} \frac{ \lambda^2}{n}+ \frac{9}{4v^4} \frac{o(\Delta_n)}{n\Delta_n}
$$
and we observe that
$$
\E_0 \left\{\I_{\{|\gamma_i|\geq \epsilon s_n\}} \gamma^2_i  \right\}
=\int\limits_{\{|\gamma_i|\geq \epsilon s_n\}} \gamma^2_n \de P_0
=
 \int\limits_{\{|\gamma_i|\geq \epsilon s_n\}} \frac{|\gamma_i|^3}{|\gamma_i|} \de P_0
\leq
\frac{1}{\epsilon s_n}\int\limits_{\{|\gamma_i|\geq \epsilon s_n\}} |\gamma_i|^3 \de P_0
\leq
\frac{1}{\epsilon s_n} \E_0 |\gamma_i|^3
$$
thus
$$
\frac{1}{s_n^2} \sum_{i=1}^n \E_0 \left\{\I_{\{|\gamma_i|\geq \epsilon s_n\}} \gamma^2_i  \right\} \leq
\frac{n}{\epsilon s_n^3}  \E |\gamma_i|^3
$$
Now we majorate $\E |\gamma_i|^3$
$$
\begin{aligned}
\E|\gamma_i|^3 =& \frac{27}{8 v^6n^3 \Delta_n^9}\E\left|(v^2\Delta_n^2 - \frac23 \lambda_0 \Delta_n^3) - \eta_i^2\right|^3\\
=&  \frac{27}{8 v^6n^3 \Delta_n^9}\left(
\left|\frac{2}{3}  \lambda_0 \Delta_n^3\right|^3 P\{N[t_{i-1},t_i)=0\} +
\left|K \Delta_n^2 + o(\Delta_n^2)\right|^3 P\{N[t_{i-1},t_i)\geq 1\} 
\right)\\
&\leq
\frac{27}{8v^6 n^3 \Delta_n^9}\left( \frac{8}{27}\lambda_0^3\Delta_n^9 + \Delta_n \left|K \Delta_n^2 + o(\Delta_n^2)\right|^3\right)\\
&\simeq
\frac{1}{n^3} + \frac{1}{n^3 \Delta_n^2}
\end{aligned}
$$
in the above we have majorated $P\{N[t_{i-1},t_i)=0\}$ by 1 and denoted by $K$ a generic constant.
Easy manipulation shows that
$$
s_n^3 = s_n^2 \sqrt{s_n^2} = \frac{1}{\sqrt{(n \Delta_n)^3}} 
$$
hence
$$
\begin{aligned}
\frac{n \E_0|\gamma_i|^3}{s_n^3} & \lesssim \frac{\Delta_n^\frac32}{\sqrt{n}} + \frac{1}{\sqrt{n \Delta_n}} \to 0
\end{aligned}
$$
with the notation $a \lesssim b$ indicating ``$a$ is majorated by something of the same order of $b$''.
Hence we have proven Lindeberg condition and asymptotic normality is established. Trivially and using exactly the same arguments, the Lindeberg condition can also be proved directly for $\sum_i \xi_i$.
\end{proof}

\section{An asymptotically efficient estimator}\label{sec4}
In the previous section we have seen that the estimators $\lambda_n^*$ and $\check\lambda_n$  are not efficient because their asymptotic variance is $\frac65 \lambda$ which is strictly greater than  $\lambda$ which is the asymptotic variance of $N(T)/T$, the asymptotically efficient estimator of the continuosly observed
 Poisson process. By a different approach we now present an asymptotically efficient estimator.
Consider the following statistic
\begin{equation}
\tilde\lambda_n = \frac{1}{n\Delta_n} \sum_{i=1}^n {\bf 1}_{\{|\eta_i|<v\Delta_n\}}=
\frac{1}{n\Delta_n} \sum_{i=1}^n {\bf 1}_{\{N([t_{i-1},t_i))\geq 1\}} 
\label{eq:bad}
\end{equation}
The statistic $\tilde\lambda_n$ is not a good estimator of $\lambda$ for fixed $\Delta_n$. Indeed,
$$
\E_0 (\tilde\lambda_n) = \frac{1}{n\Delta_n} \sum_{i=1}^n \E_0 ({\bf 1}_{\{|\eta_i|<v\Delta_n\}})=
\frac{1-e^{-\lambda_0\Delta_n}}{\Delta_n}
$$
Hence, we propose the following estimator 
\begin{equation}
\hat\lambda_n = -\frac{1}{\Delta_n} \log\left(1-\Delta_n \tilde\lambda_n\right)
\label{eq:opt}
\end{equation}
and next theorem proves that it is the efficient estimator in this context.
\begin{theorem}
Let $\Delta_n\to 0$, $n\Delta_n\to \infty$ as $n\to\infty$. Then, for all fixed $\lambda_0>0$, the estimator $\hat\lambda_n$ in \eqref{eq:opt} is consistent, asymptotically normal and attains the minimum variance,  i.e. it is asymptotically efficient
$$
\sqrt{n\Delta_n}(\hat\lambda_n - \lambda_0) \overset{d}{\longrightarrow}N(0, \lambda_0)
$$
\end{theorem}
\begin{proof}
In order to prove consistency and asymptotic normality of $\hat\lambda_n$ we first prove the same properties for $\tilde\lambda_n$.
We have seen that $\tilde\lambda_n$ in \eqref{eq:bad} is such that
$$\E_0 (\tilde\lambda_n) = \frac{1-e^{-\lambda_0\Delta_n}}{\Delta_n} = 
 \lambda_0+ \lambda_0^2 \Delta_n^2 + o(\Delta_n^2)\to \lambda_0$$
hence consistency of $\tilde\lambda_n$ follows trivially.
We now prove asymptotic normality.
Let us consider the following quantity
$$
\begin{aligned}
U_n &= \sqrt{n\Delta_n}\left(\tilde\lambda_n -\E_0 \tilde\lambda_n\right)\\
 &= \frac{1}{\sqrt{n\Delta_n}}
\sum_{i=1}^n \left\{{\bf 1}_{\{|\eta_i|<v\Delta_n\}} - \E_0\left({\bf 1}_{\{|\eta_i|<v\Delta_n\}}\right)\right\}\\
&= \frac{1}{\sqrt{n\Delta_n}}
\sum_{i=1}^n\left\{ {\bf 1}_{\{|\eta_i|<v\Delta_n\}} - \left(1-e^{-\lambda_0 \Delta_n}\right)\right\}\\
&=\sum_{i=1}^n \xi_i
\end{aligned}
$$
with
$$
\xi_i = \frac{1}{\sqrt{n\Delta_n}} \left\{{\bf 1}_{\{|\eta_i|<v\Delta_n\}} - \left(1-e^{-\lambda_0 \Delta_n}\right)\right\}
$$
We have that $\E_0 \xi_i=0$ thus $\E_0 (U_n) = 0$.
Moreover,
$$
\begin{aligned}
\V_0(\xi_i) &= \V_0 ({\bf 1}_{\{|\eta_i|<v\Delta_n\}} ) \\
&= \E_0({\bf 1}_{\{|\eta_i|<v\Delta_n\}} )(1- \E_0 ({\bf 1}_{\{|\eta_i|<v\Delta_n\}} ))\\
&= (1-e^{-\lambda_0 \Delta_n}) e^{-\lambda_0 \Delta_n} \\
&= \lambda_0 \Delta_n + o(\lambda_0 \Delta_n)
\end{aligned}
$$
hence
$$
\V_0 (U_n) = \frac{1}{n\Delta_n} n (\lambda_0 \Delta_n + o(\lambda_0\Delta_n)) = \lambda_0 + o(1)
$$
Finally, the $\xi_i$'s are independent because they only involves the absolute value of the increments $\eta_i$, hence once we prove the Lindeberg condition we also have asymptotic normality. For large $n$ it holds true that $|\xi_i|\leq 1/\sqrt{n\Delta_n}$, then Lindeberg condition is trivially true
$$
\sum_{i=1}^n \E_0 \left\{\I_{\{|\xi_i|\geq \epsilon\}} \xi^2_i  \right\} \to 0
$$
and the following result holds true
$$
U_n \overset{d}{\to} N(0,\lambda_0)
$$
Now we need to prove asymptotic normality of $\hat\lambda_n$ in \eqref{eq:opt}.
Let 
$$f(u) = -\frac{1}{\Delta_n} \log(1-u\Delta_n), \quad f'(u) = \frac{\de}{\de u} f(u) = \frac{1}{1-u\Delta_n},
$$
and
$$
\hat\lambda_n = f(\tilde\lambda_n), \quad
\lambda_0 = f\left(\tilde\lambda_0=\frac{1-e^{-\lambda_0 \Delta_n}}{\Delta_n}\right)
$$
then, by the so-called $\delta$-method,
$$
\begin{aligned}
\sqrt{n\Delta_n}(\hat\lambda_n - \lambda_0) &=
\sqrt{n\Delta_n}(f(\tilde\lambda_n) - f(\tilde\lambda_0)) \\
&=\sqrt{n\Delta_n}(\tilde\lambda_n - \tilde\lambda_0)f'(\tilde\lambda_0)+o_p(\sqrt{n\Delta_n}|\tilde\lambda_n - \tilde\lambda_0|)\\
&=\sqrt{n\Delta_n}(\tilde\lambda_n - \tilde\lambda_0) \frac{1}{1-\lambda_0\Delta_n} + o_p(1)
\end{aligned}
$$
hence
$$
\sqrt{n\Delta_n}(\hat\lambda_n - \lambda_0) \overset{d}{\to} N(0,\lambda_0)
$$
\end{proof}

\section{Conclusion}
In this paper we have shown that a consistent and asymptotically efficient estimators of the parameter $\lambda$ of the telegraph process observed at discrete times can be obtained under the high frequency sampling when the length of observation interval $[0,T]$ increases. 
We have also shown that  the moment type estimators $\lambda^*_n$ and $\check\lambda_n$ are  consistent but not efficient. They are also asymptotically Gaussian (the approximated moment type estimator $\lambda^*_n$ requires the additional assumption  $n\Delta_n^3\to 0$).
This paper also contains the explicit formula for the moments $\E X(t)^p$ of any order $p\geq 1$ of the telegraph process and a simple expansion for small $t$.

\section{Acknowledgments}
The work of the first author was supported by  JSPS (Japan Society for the Promotion of Science) Program FY2006, grant ID No. S06174.
He is also  thankful to the Graduate School of Mathematical Sciences, University of Tokyo as host research institute for the JSPS Program.

\end{document}